# Martingales and first passage times of AR(1) sequences

ALEXANDER NOVIKOV*† and NINO KORDZAKHIA ‡

†Department of Mathematical Sciences, University of Technology, Sydney, NSW 2007, Australia

‡ Department of Statistics, Macquarie University; NSW 2109, Australia

Using the martingale approach we find sufficient conditions for exponential boundedness of first passage times over a level for ergodic first order autoregressive sequences (AR(1)). Further, we prove a martingale identity to be used in obtaining explicit bounds for the expectation of first passage times.

## 1. Introduction.

We define an AR(1) sequence as a solution of the equation

$$X_n = \lambda X_{n-1} + \eta_n, \quad n = 1, 2, ... \tag{1}$$

where $\{\eta_n\}$ is a sequence of independent identically distributed random variables (*innovation*), $X_0 = x$ and $\lambda$ are nonrandom constants,

$$0 < \lambda < 1.$$

The solution of (1) has the following representation for $n = 1, 2, ...$

$$X_n = \lambda^n x + \sum_{k=0}^{n-1} \lambda^k \eta_{n-k} . \tag{2}$$

In applications such as statistical surveillance [5] and many others it is of interest to know properties of the first passage time

$$\tau_a = \inf\{n \geq 0 : X_n > a\}, \quad a \geq x,$$

where we assume $\inf\{\emptyset\} = \infty$ and so $\tau_a = \infty$ on the set $\{\sup_{t \geq 0} X_t < a\}$. Note that if

$$\eta_n \leq H \tag{3}$$

---

*Corresponding author. Email: prob@it.uts.edu.au

then in view of (2) for $n = 1, 2, ...$

$$X_n \leq \lambda^n x + H \sum_{k=0}^{n-1} \lambda^k = \lambda^n(x - \frac{H}{1-\lambda}) + \frac{H}{1-\lambda},$$

and, hence, $\tau_a = \infty$ when $H/(1-\lambda) \in (x, a)$.

In applications, the distribution and expectation of such passage times are usually approximated via Monte-Carlo simulation or using Markov chain approximations (see e.g. [16]). However, analytical bounds are also of interest (e.g. to control an accuracy of simulation algorithms).

In this paper we describe some martingales related to AR(1) sequences in the case when the innovation $\{\eta_n\}$ has a distribution with a light right tail (see Propositions 1 - 3 in Section 2). In Section 3 we use the martingale approach to show that if instead of (3) the following inequality holds

$$P\{\eta_1 > a(1-\lambda)\} > 0, \tag{4}$$

then, under some mild assumptions on the left tail of $\eta_1$, the distribution of $\tau_a$ is exponentially bounded, see Theorem 1 in Section 3. In Theorem 3, Section 4, we prove a martingale identity (the analog of classical Wald's identity for random walks) and show how to use it to obtain bounds for $E\tau_a$.

The studies related to the first passage problem for AR(1) sequences are often employed in application fields such as surveillance analysis [5], signal detection and many other areas. The martingale technique has already been used in [11] and [12] for deriving analytical approximations to the distribution and expectation of $\tau_a$ for discrete and continuous time AR(1) type processes, in particular, for the Ornstein-Uhlenbeck (O-U) processes in continuous time framework. In [13] a criterion on exponential boundedness of first passage time of O-U processes driven by general Levy processes is proposed. The papers [14] and [3] contain some analytical and numerical results on the distribution of first passage time of O-U type processes arising in engineering applications. A survey of some results in this area has been presented in lecture notes [15], where, in particular, the explicit formula for the moment generating function of first passage time of AR(1) sequences with the exponentially distributed innovation was suggested. Although, the latter case is not discussed in the current paper, we would like to mention recent papers by Jacobsen and Jensen [7], Jacobsen [8] who have obtained more general results in the case when the innovation in AR (1) process is represented in terms of a mixture of exponential random variables; furthermore, in this case the distribution of the overshoot over a level by AR(1) process was obtained.

The first passage problems for random walks and Levy processes were studied over decades. Major achievements in this direction have been obtained via the Wiener-Hopf factorisation in [2], [4],[6], [10] and many other papers and monographs. An extension of the Wiener-Hopf factorisation technique to AR(1) sequences and O-U processes is yet to be found.

**2. Some martingales associated with AR(1) sequences.**

In this paper we always consider martingales with respect to the natural filtration $\mathcal{F}_n = \sigma\{X_0, X_1, ..., X_n\}$. First, consider a martingale $M_n$ of the form

$$M_n = \lambda^{vn} q_v(X_n), \tag{5}$$

where a deterministic function $q_v(y)$ depends on a parameter $v$, the variable $y$ takes values from the domain $D$ of $X_n$. Note that, typically, $D = (-\infty, \infty)$ but if, for example, condition (3) holds and $H(1-\lambda) \geq x$ then $D \subset (-\infty, (1-\lambda)H)$.

Under the assumption that $M_n$ has a finite expectation, by definition of martingales

$$E[M_n | \mathcal{F}_{n-1}] = M_{n-1} \quad a.s.$$

which with (1) is equivalent to the equation

$$\lambda^{vn} E[q_v(\lambda X_{n-1} + \eta_n) | \mathcal{F}_{n-1}] = \lambda^{v(n-1)} q_v(X_{n-1}) \quad a.s.$$

Here $\eta_n$ is independent of $\mathcal{F}_{n-1}$ and $X_{n-1}$ may take any value from the domain $D$ of $X_n$. Therefore, if the function $q_v(y)$ is a solution of the equation

$$Eq_v(\lambda y + \eta_1) = \lambda^{-v} q_v(y), \ y \in D, \tag{6}$$

and the expectation $Eq_v(\lambda y + \eta_n)$ is finite then $\lambda^{vn} q_v(X_n)$ is a martingale.

Similar, if a martingale $M_n$ has the form

$$M_n = Q(X_n) - n \tag{7}$$

where $Q(y)$ is a deterministic function, then we obtain another equation

$$EQ(\lambda y + \eta_1) = Q(y) + 1, \quad y \in D. \tag{8}$$

Martingales of the form (5) and (7) have already been discussed in [11] and [12] under the assumptions

$$Ee^{u\eta_1} < \infty \quad for \ 0 \le u < \infty \tag{9}$$

and

$$E|\eta_1| < \infty \tag{10}$$

(though the corresponding equations (6) and (8) were not even mentioned in those papers). Here we will also use (9) but will relax condition (10) assuming only the existence of the logarithmic moment of $\eta_1$ (see Proposition 1 and 2 below) or moments of order $\delta > 0$ (see Proposition 3).

Denote the cumulant function of $\eta_1$ as follows

$$\psi(u) = \log Ee^{u\eta_1}, \ 0 \le u < \infty.$$

It is well known that if $\psi(u)$ is finite then it is a convex differentiable function for $u > 0$ (see e.g. Borovkov [2]), $\psi(0) = 0$. In view of (2) we have for any $u \in [0, \infty)$

$$Ee^{uX_n} = \exp\{\lambda^n x + \sum_{k=0}^{n-1} \psi(\lambda^k u)\}.$$

If $E\eta_n = m$ is finite then $\psi(u) = mu + o(u)$ as $u \to 0$. This fact implies that the partial sums $\sum_{k=0}^{n} \psi(\lambda^k u)$ converge to a finite limit, say, $\phi(u)$ for any $u \ge 0$ as $n \to \infty$:

$$\sum_{k=0}^{n} \psi(\lambda^k u) \to \phi(u) = \sum_{k=0}^{\infty} \psi(\lambda^k u). \tag{11}$$

Note that under assumption (9) we may have[1] $E(\eta_n^-) = \infty$ or, equivalently, $\psi'(0) = -\infty$. Under the latter condition there exists $u_0 > 0$ such that $\psi(u) < 0$ for $u \in (0, u_0)$. It implies $\psi(\lambda^k u) < 0$ for all $u \in (0, u_0)$ and for all $u \ge u_0$ and $k > \log(u/u_0)/\log(1/\lambda)$. Therefore, the series $\sum_{k=0}^{\infty} \psi(\lambda^k u)$ converges to a finite value or diverges to $-\infty$ for all $u > 0$. We now show that this series converges under the Vervaat condition

$$E \log(1 + |\eta_1|) < \infty. \tag{12}$$

---

[1] $x^- = \max(-x, 0), x^+ = \max(x, 0)$

**Lemma 1.** *Let conditions (9) and (12) hold. Then the function*

$$\phi(u) = \sum_{k=0}^{\infty} \psi(\lambda^k u)$$

*is differentiable for $u > 0$,*

$$\phi(u) = \lim_{n \to \infty} \log E e^{uX_n} = \log E e^{u\Theta}, \ \Theta \stackrel{d}{=} \sum_{k=0}^{\infty} \lambda^k \eta_{k+1}$$

*and*

$$\phi(u) = \phi(\lambda u) + \psi(u), \ 0 \leq u < \infty. \tag{13}$$

**Proof.** Accordingly to the results of Vervaat [18], under condition (12) the process $X_n$ converges in distribution as $n \to \infty$:

$$X_n \stackrel{d}{\to} \Theta, \tag{14}$$

where $\Theta$ is a finite random variable. By (2) we obtain

$$X_n = \lambda^n x + \sum_{k=0}^{n-1} \lambda^k \eta_{k+1} \stackrel{d}{\to} \Theta$$

and thus

$$\Theta \stackrel{d}{=} \sum_{k=0}^{\infty} \lambda^k \eta_{k+1}. \tag{15}$$

Since $\sum_{k=0}^{\infty} \lambda^k \eta_{k+1} = \lambda \sum_{k=1}^{\infty} \lambda^{k-1} \eta_{k+1} + \eta_1$ we have

$$\Theta \stackrel{d}{=} \lambda\Theta + \eta_1, \tag{16}$$

where in the right-hand side (RHS) random variables $\Theta$ and $\eta_1$ are assumed to be independent.

Replacing $\eta_k$ by $\eta_k^-$, we obtain from (14) and (15) as $n \to infty$

$$\sum_{k=0}^{n} \lambda^k \eta_{k+1}^- \to \sum_{k=0}^{\infty} \lambda^k \eta_{k+1}^-.$$

By the Lebesgue dominated convergence theorem this implies for any $u > 0$

$$E e^{-u \sum_{k=0}^{n} \lambda^k \eta_{k+1}^-} \to E e^{-u \sum_{k=1}^{\infty} \lambda^k \eta_k^-} > 0.$$

Since $\sum_{k=0}^{n-1} \lambda^k \eta_{n-k} = \sum_{k=1}^{n} \lambda^k \eta_k \geq -\sum_{k=1}^{n} \lambda^k \eta_k^-$ we obtain that

$$\liminf_{n \to \infty} \sum_{k=0}^{n} \psi(\lambda^k u) = \liminf_{n \to \infty} \log E e^{u \sum_{k=0}^{n} \lambda^k \eta_k} \geq \log E e^{-u \sum_{k=0}^{\infty} \lambda^k \eta_k^-} > -\infty$$

and, hence, (11) holds under condition (12).

Thus, we have shown that for any $u \in [0, \infty)$ the function

$$\phi(u) = \lim_{n \to \infty} Ee^{uX_n} = \log Ee^{u\Theta}$$

is finite. It is a differentiable function because it is a cumulant function of a finite random variable. Furthermore, in view (16) $\phi(u)$ satisfies (13). $\square$

If $E\eta_n = m$ is finite then $E\Theta = \frac{mu}{1-\lambda}$ and as $u \to 0$

$$\phi(u) = \frac{mu}{1-\lambda} + o(u). \tag{17}$$

As an illustrating example we consider the case of so-called spectrally negative stable AR(1) sequences with

$$\psi(u) = mu + Sgn(\alpha - 1)Cu^{\alpha}, \ 0 \leq u < \infty,$$

where $\alpha \in (0, 1)$ or $\alpha \in (1, 2]$, $C > 0$ (note that $C = Var(\eta_1)/2$ in the case $\alpha = 2$ which corresponds to Gaussian AR(1) sequences). One can verify that for this example

$$\phi(u) = \frac{mu}{1-\lambda} + \frac{Sgn(\alpha - 1)Cu^{\alpha}}{1-\lambda^{\alpha}}, \ 0 \leq u < \infty. \tag{18}$$

We shall often use the following condition

$$\int_1^{\infty} e^{uy - \phi(u)} u^{v-1} du < \infty \tag{19}$$

with some real $v$ (to be specified) and any $y \in D$.

A simple sufficient condition for validity of (19) is

$$\lim_{u \to \infty} \frac{\phi(u)}{u} = \infty. \tag{20}$$

Note that for the case (18) condition (12) holds for all $\alpha > 0$ and condition (19) holds when $\alpha \in (1, 2]$ but if $\alpha \in (0, 1)$ it holds only for $y \leq m$. Also, note that if the innovation $\eta_n$ is bounded from above then (20) does not hold (see Lemma 3 below) but still condition (19) could hold for $y \in D$.

Set

$$N_v(y) = \int_0^{\infty} e^{uy - \phi(u)} u^{v-1} du. \tag{21}$$

**Proposition 1.** *Let $v > 0$, conditions (9), (12) and condition (19) with $v > 0$ hold. Then*

$$\lambda^{vn} N_v(X_n) \quad \text{is a martingale.} \tag{22}$$

**Proof.** The function $N_v(y)$ is finite due the imposed conditions. Now we are going to verify that equation (6) holds for $q_v(y) = N_v(y)$. By Fubinni's theorem we have

$$EN_v(\lambda y + \eta_1) = \int_0^{\infty} \exp\{u\lambda y + \psi(u) - \phi(u)\} u^{v-1} du =$$

(with use of (13))

$$= \int_0^\infty e^{u\lambda y - \phi(\lambda u)} u^{v-1} du = \lambda^{-v} \int_0^\infty e^{uy - \phi(u)} u^{v-1} du = \lambda^{-v} N_v(y).$$

Thus, we have shown that the function $N_v(y)$ is a solution of equation (6) and hence the process $\lambda^{vn} N_v(X_n)$ is a martingale. □

Set

$$H(y) = \frac{1}{\log(1/\lambda)} \int_0^\infty (e^{uy} - 1) e^{-\phi(u)} u^{-1} du.$$

**Proposition 2.** *Let conditions (9), (12) and (19) with $v = 0$ hold. Then*

$$H(X_n) - n \quad \text{is a martingale.} \tag{23}$$

**Proof.** The function $H(y)$ is finite due the imposed conditions. By Fubinni's theorem

$$EH(\lambda y + \eta_1) = \frac{1}{\log(1/\lambda)} \int_0^\infty (e^{u\lambda y + \psi(u)} - 1) e^{-\phi(u)} u^{-1} du, \tag{24}$$

where the RHS is finite if for some $u_0 > 0$

$$\int_{u_0}^\infty |e^{u\lambda y + \psi(u)} - 1| e^{-\phi(u)} u^{-1} du < \infty \tag{25}$$

and

$$\int_0^{u_0} |e^{u\lambda y + \psi(u)} - 1| e^{-\phi(u)} u^{-1} du < \infty. \tag{26}$$

The integral $\int_{u_0}^\infty$ in (25) is finite for any $u_0 > 0$ due (19) with $v = 0$; the integral $\int_0^{u_0}$ in (26) is finite if and only if

$$\int_0^{u_0} |1 - e^{\psi(u)}| u^{-1} du < \infty. \tag{27}$$

If $E\eta_n = m$ is finite then in view of (17), obviously, condition (27) holds. If $E\eta_n^- = \infty$ then there exists $u_0 > 0$ such that $\phi(u) \leq 0$ for $u \in [0, u_0]$. This fact together with Lemma 1 lead to the following estimates

$$\int_0^{u_0} |1 - e^{\psi(u)}| u^{-1} du = \int_0^{u_0} (1 - e^{\psi(u)}) u^{-1} du =$$

$$E \int_0^{u_0} (1 - e^{u\eta_1}) u^{-1} du \leq EI\{\eta_1 \leq 0\} \int_0^{-\eta_1 u_0} (1 - e^{-u}) u^{-1} du$$

$$\leq P\{-1 \leq \eta_1 \leq 0\} \int_0^{u_0} (1 - e^{-u}) u^{-1} du +$$

$$P\{\eta_1 < -1\}(\int_0^{u_0}(1-e^{-u})u^{-1}du + EI\{\eta_1 < -1\}\int_{u_0}^{-\eta_1 u_0} u^{-1}du)$$

$$\leq const + E\log(1+|\eta_1|).$$

Thus, we have shown that integrals in (25) and (26) are finite and therefore the RHS of (24) can be now written (with use of (13)) as follows:

$$EH(\lambda y + \eta_n) = H(y) + \int_0^\infty \frac{(e^{u\lambda y - \phi(\lambda u)} - e^{uy-\phi(u)})}{\log(1/\lambda)} u^{-1} du.$$

To satisfy equation (8) we need to show that the last integral equals to 1. In fact, this type of integrals is well known; it is called the Frullani's integral (see e.g. [9]). It equals really to 1 if it is absolutely convergent and the function $\phi(u)$ is continuous (the latter is, of course, true). The absolute convergence can be checked similarly to the verifications of (25) and (26) as it had been done above. □

For the proof of Theorem 1 below we shall need also exponential martingales like (22) but with negative $v < 0$ and a modified martingale function.

Set

$$W_v(x) := \int_0^\infty (e^{ux-\phi(u)} - 1)u^{v-1}du. \tag{28}$$

**Proposition 3.** *Let conditions (9) and (19) with $v < 0$ hold. If there exists $\delta \in (0,1]$ such that*

$$E(\eta_1^-)^\delta < \infty \tag{29}$$

*then for $v \in (-\delta, 0)$*

$$\lambda^{vn} W_v(X_n) \text{ is a martingale.}$$

**Proof.** First we note that applying for the random variable $\Theta$ from Lemma 1 the following inequality

$$(\sum_{k=0}^\infty |x_k|)^\delta \leq \sum_{k=0}^\infty |x_k|^\delta, \quad \delta \in (0,1],$$

we obtain

$$E|\Theta|^\delta \leq \sum_{k=0}^\infty \lambda^{k\delta} E|\eta_k|^\delta = \frac{E|\eta_k|^\delta}{1-\lambda^\delta}.$$

Hence, $E|\Theta|^\delta$ is finite under condition (29).

In view of (19) the function $W_v(x)$ is finite if there exists $u_0 > 0$ such that

$$\int_0^{u_0} |1 - e^{\phi(u)}| u^{v-1} du < \infty. \tag{30}$$

If $\delta = 1$ (that is when $E\eta_n = m$ is finite) then in view of (17) condition (30) holds for any $u_0 > 0$ and any $v \in (-1, 0]$. If $E\eta_n^- = \infty$ then $E\Theta^- = \infty$ and due to convexity of the function $\phi(u)$ there exists $u_0 > 0$ such that $\phi(u) \leq 0$ for $u \in [0, u_0]$. This fact together with Lemma 1 leads to the following estimate

$$\int_0^{u_0} |1 - e^{\phi(u)}|u^{v-1}du = E\int_0^{u_0} (1 - e^{u\Theta})u^{v-1}du$$

$$\leq E\int_0^{u_0} (1 - e^{-u\Theta^-})u^{v-1}du \leq E(\Theta^-)^{-v}\int_0^\infty (1 - e^{-u})u^{v-1}du.$$

which is finite under condition (29) for $v \in [-\delta, 0]$.

Now by Fubinni's theorem we have

$$EW_v(\lambda y + \eta_n) = \int_0^\infty (e^{u\lambda y - \phi(\lambda u)} - 1)u^{v-1}du$$

$$= \lambda^{-v}\int_0^\infty (e^{uy - \phi(u)} - 1)u^{v-1}du = \lambda^{-v}W_v(y)$$

and, hence, $W_v(x)$ satisfies the martingale equation (6).□

**3. Exponential boundedness of first passage times.**

**Theorem 1.** *Let $a > 0$, conditions (4) and (29) with some $\delta \in (0, 1)$ hold. Then there exists $\alpha > 0$ such that*

$$Ee^{\alpha\tau_a} < \infty.$$

The main idea of the proof consists in use of a martingale $\lambda^{vn}W_v(X_n)$ with $W_v(y)$ from (28) and the optional stopping theorem to derive the following bound for some $v < 0$ and all $n \geq 1$

$$E\lambda^{v\min(\tau_a, n)} \leq C < \infty, \qquad (31)$$

where $C$ is a constant (not depending on $n$). The latter estimate along with Fatou's lemma implies Theorem 1.

To derive (31), first we note that if we truncate the positive part of the innovation $\eta_n$ from above by a constant, say, $N > 0$ then in (2) the corresponding stopping time $\tau_a$ will be greater than the original one (without truncation). So, it is sufficient to prove (31) only for the truncated innovation.

The derivation of (31) requires the following estimate for the corresponding cumulant function $\phi(u)$ associated with the truncated innovation $\{-\eta_n^- + NI\{\eta_n \geq N\}\}$.

**Lemma 3.** *Let condition (12) hold,*

$$P\{\eta_1 = N\} = p = 1 - P\{\eta_1 \leq 0\} > 0. \qquad (32)$$

*Then as $u \to \infty$*

$$\phi(u) = \frac{uN}{1 - \lambda} + o(u).$$

**Proof.** We have

$$\psi(u) = \log(Ee^{u\eta_1}) = uN - g(u) \qquad (33)$$

with

$$g(u) = -\log[EI\{\eta_1 \leq 0\}\exp\{u(\eta_1 - N)\} + p] \geq 0.$$

Note

$$g'(u) = -\frac{EI\{\eta_1 \leq 0\}(\eta_1 - N)\exp\{u(\eta_1 - N)\}}{EI\{\eta_1 \leq 0\}\exp\{u(\eta_1 - N)\} + p} \to 0 \ \ as \ \ u \to \infty.$$

In view of (33) we have

$$\phi(u) = \sum_{k=0}^{\infty} \psi(\lambda^k u) = \frac{u\,N}{1-\lambda} - \Delta(u)\,, \tag{34}$$

where

$$\Delta(u) := \sum_{k=0}^{\infty} g(\lambda^k u) \geq 0,\ 0 \leq u < \infty.$$

The function $\Delta(u)$ is differentiable and concave since by Lemma 1 the function $\phi(u)$ is differentiable and convex. Obviously,

$$\Delta(u) = \Delta(\lambda u) + g(u),\ 0 \leq u < \infty. \tag{35}$$

From concavity of $\Delta(u)$ it follows that the derivative $\Delta'(u)$ is non-increasing. It certainly has a lower bound (as $\Delta(u) \geq 0$) and therefore there exists a finite limit

$$\lim_{u \to \infty} \Delta'(u) = A.$$

Applying L'Hospitale's rule we have

$$A = \lim_{u \to \infty} \frac{\Delta(u)}{u} = \lim_{u \to \infty} \frac{\Delta(\lambda u) + g(u)}{u} = \lambda A + \lim_{u \to \infty} \frac{g(u)}{u}.$$

As $\lim_{u \to \infty} \frac{g(u)}{u} = \lim_{u \to \infty} g'(u) = 0$, we obtain that $A = \lambda A$ and therefore $A = 0$.

Lemma 3 is proved. $\square$

**Proof of Theorem 1.** Due to condition (4) we can always chose $N > a(1-\lambda) > 0$ such that

$$P\{\eta_1 > N\} > 0.$$

Now let $\tau_a$ is the first passage time for AR(1) processes $X_n$ generated by the truncated innovation with the property (32) and therefore the corresponding cumulant function has the property (32).

Note [1] that

$$X_{n \wedge \tau_a} = \lambda X_{n \wedge \tau_a - 1} + \eta_{n \wedge \tau_a} \leq \lambda a + N < \frac{N}{1-\lambda},$$

---

[1] $n \wedge \tau_a = \min(n, \tau_a)$

Now one can check that conditions of Proposition 3 hold and so the process $\lambda^{-vn}W_v(X_n)$ is a martingale. By the optional stopping theorem for any $n \geq 1$

$$E(\lambda^{v(n\wedge\tau_a)}W_v(X_{n\wedge\tau_a})) = W_v(x).$$

Since $X_{n\wedge\tau_a} \leq \lambda a + N$ and $W_v(x)$ is an increasing function of $x$ we have for any $v \in (-\delta, 0)$

$$E(\lambda^{v(n\wedge\tau_a)}W_v(\lambda a + N)) \geq W_v(x). \tag{36}$$

Note

$$W_v(x) = \int_0^\infty (e^{ux-\phi(u)} - 1)u^{v-1}du = C(x,v) + \frac{1}{v} \tag{37}$$

where

$$C(x,v) = \int_0^1 (e^{ux-\phi(u)} - 1)u^{v-1}du + \int_1^\infty e^{ux-\phi(u)}u^{v-1}du.$$

It is easy to check that $C(x,\mu) \to C(x,0)$ as $v \to 0$ and $C(x,0)$ is finite for any $x < \frac{N}{1-\lambda}$. It implies that there exists $v \in (-\delta, 0)$

$$-2|C(x,0)| + \frac{1}{v} \leq W_\mu(x) \leq W_\mu(\lambda a + N) \leq 2|C(\lambda a + N, 0)| + \frac{1}{v}$$

and thus with (36) we obtain

$$E(\lambda^{-vn\wedge\tau_a}(1 + 2v|C(\lambda a + N, 0)|)) \leq 1 - 2v|C(x,0)|.$$

Now choose sufficiently small $v < 0$ such that $1 + 2v|C(\lambda a + N, 0)| > 0$. Then this implies that for any $n \geq 1$

$$E(\lambda^{vn\wedge\tau_a}) \leq \frac{1 - 2v|C(x,0)|}{1 + 2v|C(\lambda a + N, 0)|} < \infty.$$

The proof is completed.

### 4. Bounds for the expectation of first passage times.

Theorem 1 implies that $E(\tau_a)$ is finite under conditions (4) and (29). Actually, the latter condition can be slightly relaxed.

**Theorem 2.** *Let $a > 0$, conditions (12) and (4) hold. Then*

$$E(\tau_a) < \infty.$$

**Proof.** We use Lemma 2 and the fact that the first passage time $\tilde\tau_a$ of AR(1) with a truncated innovation is greater than the original one (without truncation).

As in the proof of Theorem 1 let $N > a(1 - \lambda) > 0$ such that $P\{\eta_1 > N\} > 0$ and let $\tilde\tau_a$ is the stopping stopping time for AR(1) processes generated by the innovation with the property

$$P\{\tilde\eta_1 = N\} = p = 1 - P\{\tilde\eta_1 \leq 0\} > 0.$$

Then by Lemma 2 the corresponding cumulant function $\phi(u) = \frac{uN}{1-\lambda} + o(u)$ as $u \to \infty$.

By Proposition 2 (or Theorem 2)

$$E\tilde{\tau}_a \wedge n = \frac{1}{\log(1/\lambda)} \int_0^\infty (Ee^{u\tilde{X}_{\tilde{\tau}_a \wedge n}} - e^{ux})e^{-\phi(u)}u^{-1}du$$

where $\tilde{X}_{\tilde{\tau}_a \wedge n} \leq \lambda a + N$. It implies

$$E\tilde{\tau}_a \wedge n \leq \frac{1}{\log(1/\lambda)} \int_0^\infty (e^{u(\lambda a+N)} - e^{ux})e^{-\frac{N}{1-\lambda}(u)+o(u)}u^{-1}du = const < \infty.$$

By Fatou's lemma it implies that $E\tilde{\tau}_a < \infty$ and thus the proof is completed.

Now we prove a general martingale identity which can be used for derivation of bounds and asymptotics for $E(\tau_a)$.

**Theorem 3.** *Let conditions (12) and (19) with $v = 0$ hold. If $E(\tau_a) < \infty$ then*

$$E\tau_a = \frac{1}{\log(1/\lambda)} \int_0^\infty (Ee^{uX_{\tau_a}} - e^{ux})e^{-\phi(u)}u^{-1}du.$$

**Proof.** By Proposition 2 and the optional stopping theorem we have for any $n = 1, 2, ...$

$$E(\tau_a \wedge n) = EH(X_{\tau_a \wedge n}) - H(x) = \frac{1}{\log(1/\lambda)} E \int_0^\infty (e^{u\,X_{\tau_a \wedge n}} - e^{ux})e^{-\phi(u)}u^{-1}du\ .$$

Since $\lim_{n\to\infty} E(\tau_a \wedge n) = E(\tau_a)$ and

$$\int_0^\infty (Ee^{uX_{\tau_a \wedge n}} - e^{ux})e^{-\phi(u)}u^{-1}du$$

$$= EI(\tau_a \leq n) \int_0^\infty (e^{uX_{\tau_a}} - e^{ux})e^{-\phi(u)}u^{-1}du$$

$$+ EI(\tau_a > n) \int_0^\infty (e^{uX_n} - e^{ux})e^{-\phi(u)}u^{-1}du\ , \tag{38}$$

where the first term in RHS is a monotonic function of $n$. Therefore, we need only to show a convergence as $n \to \infty$ to zero for the latter integral term. Note that $X_n \leq a$ on the set $\{\tau_a > n\}$ and this implies the following upper bound

$$EI(\tau_a > n) \int_0^\infty (e^{uX_n} - e^{ux})e^{-\phi(u)}u^{-1}du$$

$$\leq P(\tau_a > n) \int_0^\infty (e^{ua} - e^{ux})e^{-\phi(u)}u^{-1}du \to 0$$

because $P(\tau_a > n) \to 0$ and the last integral is finite due to the imposed conditions.

To show that the lower bound for the second integral term in the RHS of (38) tends also to zero as $n \to \infty$, we note

$$EI(\tau_a > n) \int_0^\infty (e^{ux} - e^{uX_n})e^{-\phi(u)}u^{-1}du \leq$$

$$P(\tau_a > n) \int_0^\infty (e^{ux} - 1) e^{-\phi(u)} u^{-1} du + EI(\tau_a > n) \int_0^\infty (1 - e^{-uX_n^-}) e^{-\phi(u)} u^{-1} du.$$

The first integral in the RHS converges, obviously, to zero.

Note that in view of (2)

$$X_n^- \leq \lambda^n x^- + \sum_{k=0}^{n-1} \lambda^k \eta_{n-k}^- := Z_n, \quad n = 1, 2, .$$

where $Z_n$ is a AR(1) process with innovation sequence $\{\eta_k^-\}, Z_0 = x^-$.

Since

$$\int_0^\infty (1 - e^{-uX_n^-}) e^{-\phi(u)} u^{-1} du \leq \int_1^\infty e^{-\phi(u)} u^{-1} du + \int_0^1 (1 - e^{-uZ_n}) e^{-\phi(u)} u^{-1} du,$$

(where $\int_1^\infty e^{-\phi(u)} u^{-1} du$ is finite by condition (19) with $v = 0$), we need only to verify that

$$\lim_{n \to \infty} EI(\tau_a > n) \int_0^1 (1 - e^{-uZ_n}) u^{-1} du = 0.$$

Note

$$EI(\tau_a > n) \int_0^1 (1 - e^{-uZ_n}) u^{-1} du = EI(\tau_a > n) \int_0^{Z_n} (1 - e^{-u}) u^{-1} du$$

$$\leq P(\tau_a > n, Z_n \leq 1) \int_0^1 (1 - e^{-u}) u^{-1} du + EI(\tau_a > n, Z_n > 1) \int_1^{Z_n+1} (1 - e^{-u}) u^{-1} du.$$

The first term in the RHS here tends, obviously, to zero. For the second term we have

$$EI(\tau_a > n) \int_1^{Z_n+1} (1 - e^{-u}) u^{-1} du \leq EI(\tau_a > n) \log(Z_n + 1).$$

Since

$$\log(Z_n + 1) \leq \log(1 + x^- + \sum_{k=1}^n \eta_k^-) \leq \log(1 + x^-) + \sum_{k=1}^n \log(1 + \eta_k^-)$$

and by the Wald identity $E(\sum_{k=1}^{\tau_a} \log(1 + \eta_k^-)) = E(\log(1 + \eta_1^-)) E\tau_a < \infty$, we obtain due to the Lebesgue dominated convergence theorem that as $n \to \infty$

$$EI(\tau_a > n) \log(Z_n + 1) \leq EI(\tau_a > n) \log(Z_{\tau_a} + 1) \to 0$$

Now combining all estimates obtained above we complete the proof.

**Remark.** Theorem 3 can be used for obtaining bounds and asymptotics of $E\tau_a$. Since the overshoot

$$\xi_a = X_{\tau_a} - a$$

is always nonnegative, under the assumption of Theorem 3 we obtain the following lower bound

$$E\tau_a \geq \frac{1}{\log(1/\lambda)} \int_0^\infty (e^{ua} - e^{ux})e^{-\phi(u)}u^{-1}du.$$

The upper bound can be obtained with use of truncation (if the original innovation is not bounded) of the innovation sequence $\{\eta_k\}$ from above by a constant, say, $H$. For the latter case, noting that $\tilde{X}_{\tau_a} \leq \lambda a + H$, we obtain that

$$E\tau_a \leq E\tilde{\tau}_a \leq \frac{1}{\log(1/\lambda)} \int_0^\infty (e^{u(\lambda a + H)} - e^{ux})e^{-\tilde{\phi}(u)}u^{-1}du$$

where the function $\tilde{\phi}(u)$ is the corresponding cumulant of the limit distribution of the AR(1) sequence $\tilde{X}_n$ with the truncated innovation.

To complete the exposition we just note that some applications of these bounds for finding optimal designs in statistical quality control have been presented in [17].